\documentclass{elsart}
\usepackage{amssymb,amsfonts}

\newcommand{\starbar}{
\bigskip

\begin{center}
$*\qquad * \qquad *$
\end{center}

\bigskip

}

\newcommand{\trace}{\mathop{\rm Tr}\nolimits}

\newcommand{\N}{{\mathbb{N}}}
\DeclareRobustCommand\openone{\leavevmode\hbox{\small1\normalsize\kern-.33em1}}
\newcommand{\identity}{\mathrm{\openone}}
\newcommand{\id}{\identity}

\newcommand{\be}{\begin{equation}}
\newcommand{\ee}{\end{equation}}
\newcommand{\bea}{\begin{eqnarray}}
\newcommand{\eea}{\end{eqnarray}}
\newcommand{\beas}{\begin{eqnarray*}}
\newcommand{\eeas}{\end{eqnarray*}}

\newtheorem{theorem}{Theorem}
\newtheorem{lemma}{Lemma}

\begin{document}
\begin{frontmatter}
\title{A Lieb-Thirring inequality for singular values}
\author{Koenraad M.R. Audenaert}
\address{
Mathematics Department\\
Royal Holloway, University of London\\
Egham TW20 0EX, United Kingdom
}
\ead{koenraad.audenaert@rhul.ac.uk}
\date{\today}
%------------------------------------------------------------------ ABSTRACT
\begin{abstract}
Let $A$ and $B$ be positive semidefinite matrices.
We investigate the conditions under which the Lieb-Thirring inequality can be extended to singular values.
That is, for which values of $p$ does the majorisation $\sigma(B^pA^p) \prec_w \sigma((BA)^p)$ hold, 
and for which values its reversed inequality $\sigma(B^pA^p) \succ_w \sigma((BA)^p)$.
\end{abstract}

\end{frontmatter}
%--------------------------------------------------------------------------------------------------
The famous Lieb-Thirring inequality \cite{lieb} states that for positive semidefinite matrices $A$ and $B$,
and $p\ge1$, $\trace(AB)^p\le \trace(A^p B^p)$, while for $0<p\le 1$ the inequality is reversed.
Many generalisations of this inequality exist \cite{ka,wang}, 
one of the most noteable being the Araki-Lieb-Thirring
inequality \cite{araki}. For positive matrices $A$ and $B$, and any unitarily invariant norm $|||\cdot|||$, 
the following holds (see also Theorem IX.2.10 in \cite{bhatia}):
$|||(BAB)^p||| \le |||B^p A^p B^p|||$ when $p\ge 1$, and the reversed inequality when $0<p\le 1$.
This inequality can be equivalently expressed as the weak majorisation relation between singular values
$\sigma((BAB)^p) \prec_w \sigma(B^p A^p B^p)$. Here, $\sigma(X) \prec_w \sigma(Y)$ if and only if
$\sum_{j=1}^k \sigma_j(X) \le \sum_{j=1}^k \sigma_j(Y)$, for $1\le k\le d$, 
where $\sigma_j(X)$ denotes the $j$-th largest singular value of $X$.

In this paper we study the related question whether a majorisation relation exists between the singular values
of the non-symmetric products $B^p A^p$ and $(BA)^p$. 
The latter expression is well-defined because the eigenvalues of a product
of positive semidefinite matrices are real and non-negative.
Our main result is the following:
\begin{theorem}\label{th:majo0}
Let $A,B\ge0$ be $d\times d$ matrices. For $0<p\le 1/2$,
\be
\sigma(B^pA^p) \prec_w \sigma((BA)^p). \label{eq:ltsa}
\ee
In addition, if $d=2$, the range of validity extends to $0<p\le 1$.

For $p\ge d-1$ and for $p\in\N_0$, the reversed inequality holds:
\be
\sigma(B^pA^p) \succ_w \sigma((BA)^p). \label{eq:ltsb}
\ee
\end{theorem}

In the first half of the paper, we prove this Theorem for $p$ satisfying the condition 
$1/p\in\N_0$ or $1/p\ge d-1$ or $p\in\N_0$ or $p\ge d-1$.
We do so by chaining together two majorisations; in terms of the first inequality, (\ref{eq:ltsa}), we chain together
$\sigma(A^p B^p)\prec_w\sigma^p(AB)$ and $\sigma^p(AB) \prec_w \sigma((AB)^p)$.
While the first majorisation indeed holds generally and is a straightforward consequence of the original
Lieb-Thirring inequality, see Theorem \ref{th:majo1},
the second majorisation turns out to be subject to the rather surprising condition on $p$ (Theorem \ref{th:majo2}).
In the second half of this paper, we follow a different route and obtain validity of (\ref{eq:ltsa}) 
for $0<p\le 1/2$.

Henceforth, we abbreviate the term positive semidefinite as PSD.

The following Theorem is already well-known:
\begin{theorem}\label{th:majo1}
For $A,B$ PSD, and $0<p\le 1$,
$$
\sigma(A^p B^p)\prec_w\sigma^p(AB).
$$
For $p\ge1$, the direction of the majorisation is reversed.
\end{theorem}
\textit{Proof.}
We only have to prove the statement for $\sigma_1$, i.e.\ the infinity norm $||.||$.
From that we can derive the full majorisation statement
by using the well-known trick, due to Weyl, of replacing $X$ by its antisymmetric tensor powers, 
as in \cite{araki}.

Consider $0<p\le1$.
By the Araki-Lieb-Thirring inequality for the infinity norm $||\cdot||$, we have
$$
||AB^2A||^p = ||(AB^2A)^p|| \ge ||A^pB^{2p}A^p||.
$$
Noting that $||XX^*||=||X||^2$, this gives, indeed,
$$
||AB||^p \ge ||A^p B^p||.
$$
This inequality was first proven by Heinz (see Theorem IX.2.3 in \cite{bhatia}).
For $p\ge1$, the direction of the inequalities is reversed.
\qed

For the second majorisation we need a lemma, which relates the question to a result by FitzGerald and Horn.
\begin{lemma}\label{lem:hfg}
Let $(\lambda_i)_i$ be a sequence of $d$ non-negative numbers.
The $d\times d$ matrix $C$ with entries
$$
C_{i,j} = \frac{1 - \lambda_i^{\alpha}\lambda_j^{\alpha}}{1-\lambda_i\lambda_j},
$$
is PSD if $\alpha\in \N_0$ or $\alpha\ge d-1$.
\end{lemma}
\textit{Proof.}
This expression can be represented in integral form as \cite{HJII}
$$
C_{i,j} = \alpha \int_0^1 dt\,(t+(1-t)\lambda_i\lambda_j)^{\alpha-1}.
$$
Thus $C$ is PSD if the integrand is.
Since for $0\le t\le 1$ the matrix $(t+(1-t)\lambda_i\lambda_j)_{i,j}$ is PSD and has non-negative entries,
$C$ being PSD follows from a Theorem of FitzGerald and Horn \cite{HJII} that states that the $q$-th entrywise
power of an entrywise non-negative PSD matrix is again PSD, provided either $q\in\N_0$ or $q\ge d-2$.
Here, $q=\alpha-1$, hence the condition is $\alpha\in\N_0$ or $\alpha\ge d-1$.
\qed

\begin{theorem} \label{th:majo2}
Let $X$ be a $d\times d$ matrix with non-negative real eigenvalues.
For $p$ in the range $0<p\le1$, the majorisation
$$
\sigma^p(X) \prec_w \sigma(X^p)
$$
holds, provided $1/p\in\N_0$ or $1/p\ge d-1$.

For the range $p\ge 1$, the direction of the majorisation is reversed, and the conditions for validity are
$p\in\N_0$ or $p\ge d-1$.
\end{theorem}
\textit{Proof.}
Consider the case $0<p\le 1$ first.

Again, we consider the inequality $\sigma_1^p(X) \le \sigma_1(X^p)$,
from which the majorisation of the Theorem follows
by the Weyl trick.

An equivalent statement of the inequality is:
$||X^p||=1$ implies $||X||\le 1$ (obtainable via rescaling $X$).

If we impose that $X$ be diagonalisable, it has an eigenvalue decomposition $X=S\Lambda S^{-1}$,
where $S$ is invertible, and $\Lambda$ is diagonal, with diagonal entries $\lambda_k\ge0$.
Then
\beas
||X^p||=1 &\Longrightarrow& (X^p)^* (X^p)\le\id \\
&\Longrightarrow& S^{-*} \Lambda^p S^*S\Lambda^p S^{-1}\le\id \\
&\Longrightarrow& \Lambda^p S^*S\Lambda^p \le S^*S.
\eeas
Let us introduce the matrix $A=S^*S$, which of course is positive definite, by invertibility of $S$.
Thus the statement $||X^p||=1$ is equivalent with $\Lambda^p A \Lambda^p\le A$.
Likewise, the statement $||X||=1$ is equivalent with $\Lambda A \Lambda\le A$.
We therefore have to prove the implication
\be
\Lambda^p A \Lambda^p\le A \Longrightarrow \Lambda A \Lambda\le A.
\label{eq:imp1}
\ee
%Note: A direct consequence of $\Lambda^p A \Lambda^p\le A$ is $\Lambda\le\id$, since the diagonal
%entries of $\Lambda^p A \Lambda^p$ are $\lambda_k^{2p} A_{kk}$. By the inequality, these have to be
%less or equal $A_{kk}$, whence $\lambda_k^{2p}$, and therefore $\lambda_k$, has to be less or equal 1.

Now note that, since $\Lambda$ is diagonal, the condition $\Lambda^p A \Lambda^p\le A$ can be written
as
$$
A':=A\circ (1-\lambda_i^p\lambda_j^p)_{i,j=1}^d \ge 0,
$$
where $\circ$ denotes the Hadamard product.
Likewise, $\Lambda A \Lambda\le A$ can be written as
$$
A\circ (1-\lambda_i\lambda_j)_{i,j=1}^d \ge 0.
$$
In terms of $A'$, this reads
$$
A'\circ C \ge 0,
$$
with
$$
C:=\left(\frac{1-\lambda_i\lambda_j}{1-\lambda_i^p\lambda_j^p}\right)_{i,j=1}^d.
$$
Thus, by Schur's Theorem \cite{HJI},
the implication (\ref{eq:imp1}) would follow from non-negativity of the matrix $C$.
Using Lemma \ref{lem:hfg}, we find that a sufficient condition is $1/p\in\N_0$ or $1/p\ge d-1$.

Using a standard continuity argument, we can now remove the restriction that $X$ be diagonalisable.

The case $p>1$ is treated in a completely similar way, but relying instead on the non-negativity of the matrix 
$$
\left(\frac{1-\lambda_i^p\lambda_j^p}{1-\lambda_i\lambda_j}\right)_{i,j=1}^d.
$$
\qed

For all other values of $p$ than the mentioned ones, the matrix $C$ encountered in the proof 
is in general no longer non-negative.
Likewise, for these other values of $p$, counterexamples can be found to the inequality that we wanted to prove here,
so the given conditions on $p$ are the best possible.

Combining Theorem 2 and Theorem 3 immediately proves Theorem 1 for $1/p\in\N_0$ or $1/p\ge d-1$
or $p\in\N_0$ or $p\ge d-1$.

\starbar %%%%%%%%%%%%%%%%%%%%%%%%%%%%%%%%%%%%%%%%%%%%%%%%%%%%%%%%%%%%%%%%%%%%%%%%%%%%%

To prove the remaining case covered by Theorem \ref{th:majo0}, 
we derive several equivalent forms of the inequalities (\ref{eq:ltsa}) and (\ref{eq:ltsb}).
We again
only need to treat the $\sigma_1$ case, as the full
statement follows from it using the Weyl trick. 

Consider first the case $0< p\le 1$. Then we need to consider
\be
|| B^pA^p || \le ||(BA)^p||,
\label{eq:baba}
\ee
since the largest singular value is just the operator norm.

As a first step, we reduce the expressions in such a way that only positive matrices appear with a fractional power.

By exploiting the relation $||X||=||X^* X||^{1/2}$, (\ref{eq:baba}) is equivalent to
$$
|| A^p B^{2p} A^p || \le ||(AB)^p(BA)^p||,
$$
which, by homogeneity of both sides, can be reformulated as
$$
||(AB)^p(BA)^p|| \le 1 \Longrightarrow || A^p B^{2p} A^p || \le 1,
$$
and, in terms of the PSD ordering,
\be\label{eq:abba1}
(AB)^p(BA)^p \le \id \Longrightarrow A^p B^{2p} A^p \le \id.
\ee

\begin{lemma}
For any $A>0$ and $B\ge0$, there exist diagonal $\Lambda\ge0$ and invertible $S$ such that
$A=SS^*$ and $AB=S\Lambda S^{-1}$, and, consequently, $B=S^{-*}\Lambda S^{-1}$.
\end{lemma}
\textit{Proof.}
Let $AB=T\Lambda T^{-1}$ be an eigenvalue decomposition of $AB$. Because $A$ and $B$ are PSD, the eigenvalues of $AB$ are non-negative,
hence $\Lambda\ge0$. Assuming that all eigenvalues of $AB$ are distinct, we show that $T^{-1}AT^{-*}$ is necessarily diagonal.

Indeed, from $AB=T\Lambda T^{-1}$ follows $T^{-1}AT^{-*}\,\,T^*BT=\Lambda$. 
The factors $X=T^{-1}AT^{-*}$ and $Y=T^*BT$ are positive definite,
and positive semidefinite, respectively, since they are related to $A$ and $B$ by a $*$-conjugation.
Now note that $\Lambda$ is diagonal and all its diagonal elements are distinct.
This implies that $X$ and $Y$, both Hermitian, are themselves diagonal.
This follows from taking the hermitian conjugate of $XY=\Lambda$, $YX=\Lambda$, and noting that the two equations
taken together imply that $X$ and $Y$ commute and are therefore diagonalised by the same unitary conjugation.
Then we see that the product $XY$ must also be diagonalised by that same unitary conjugation. However, $XY=\Lambda$ is already diagonal,
so that $X$ and $Y$ must be diagonal too.

By a continuity argument, we see that there must exist a $T$ diagonalising both $AB$ (via a similarity) and $A$ (via a $*$-conjugation)
even when the eigenvalues of $AB$ are not distinct.

The lemma now follows by putting $S=T X^{1/2}$.
\qed

Using the Lemma, the left-hand side (lhs) of (\ref{eq:abba1}) can be rewritten as
\beas
(AB)^p(BA)^p &=& (S\Lambda S^{-1})^p(S^{-*}\Lambda S^*)^p \\
&=& S\Lambda^p S^{-1} S^{-*} \Lambda^p S^*.
\eeas
The condition $(AB)^p(BA)^p\le \id$ then becomes
$$
\Lambda^p S^{-1} S^{-*} \Lambda^p \le S^{-1} S^{-*},
$$
which turns into
$$
\Lambda^p C \Lambda^p \le C
$$
on defining $C=S^{-1} S^{-*}>0$.

Similarly, the condition of the right-hand side (rhs) of (\ref{eq:abba1}), 
$A^p B^{2p} A^p\le\id$, can be rewritten as $B^{2p} \le A^{-2p}$, or
\be
(S^{-*}\Lambda S^{-1})^{2p} \le (SS^*)^{-2p} = (S^{-*}S^{-1})^{2p}.
\ee
Using the polar decomposition, we can put $S^{-*} = UC^{1/2}$, where $U$ is a unitary matrix.
Then the condition of the rhs becomes $(UC^{1/2} \Lambda C^{1/2}U^*)^{2p} \le (UCU^*)^{2p}$, or
\be
(C^{1/2} \Lambda C^{1/2})^{2p} \le C^{2p}.
\ee

Thus, implication (\ref{eq:abba1}) is equivalent to
\be\label{eq:lcl1}
\Lambda^p C \Lambda^p \le C \Longrightarrow
(C^{1/2} \Lambda C^{1/2})^{2p} \le C^{2p},
\ee
for $0\le p\le 1$, and $C>0$, $\Lambda\ge0$.

On left- and right-multiplying both sides of the lhs of (\ref{eq:lcl1}) with $C^{1/2}$,
we get
$$
(C^{1/2} \Lambda^p C^{1/2})^2 \le C^2 \Longrightarrow (C^{1/2} \Lambda C^{1/2})^{2p} \le C^{2p}.
$$
By putting $A=C^{1/2}$ and $B=\Lambda^{p}$, this becomes
$$
(ABA)^2 \le A^4 \Longrightarrow (A B^{1/p} A)^{2p} \le A^{4p}.
$$
In this equivalent form, it is now easy to prove (\ref{eq:ltsa}) for $p\le 1/2$.

\textit{Proof of Theorem 1 for $0\le p\le1/2$:}
By operator monotonicity of the square root, $(ABA)^2\le A^4$ implies $ABA\le A^2$.
Dividing out $A$ on both sides, this is equivalent with $B\le \id$.
This implies $B^{1/p}\le \id$, for all $p>0$, and thus $AB^{1/p}A\le A^2$.
Since $0<p\le 1/2$, operator monotonicity of the $2p$-th power
finally implies $(AB^{1/p}A)^{2p}\le A^{4p}$.
\qed

\medskip

For $d>2$ and $1/2<p<1$, we have found counterexamples.
To narrow down the search for counterexamples, we semi-intelligently 
chose a random positive diagonal $d\times d$
matrix $D$ and a random $d$-dimensional vector $\psi$ to construct $A$ and $B$ matrices:
\beas
A^2 &=& \left(\frac{\psi_k\overline{\psi}_l}{1-D_{kk}D_{ll}}\right)_{i,j=1}^d \\
B &=& ||A^{-1}DA^2D A^{-1}||^{-1/2}\,D.
\eeas
The condition $(ABA)^2 \le A^4$ is equivalent with $||A^{-1}BA^2BA^{-1}||\le 1$
and is thus satisfied by construction.
However, with high probability $A$ and $B$ are found that violate $(A B^{1/p} A)^{2p} \le A^{4p}$.
As the violations are extremely small, all calculations have to be done in high-precision arithmetic
(we used 60 digits of precision)
\footnote{A Mathematica notebook with these calculations is available from the author on request.}.
This numerical procedure yielded counterexamples for $d=3$ and $p$ between 0.89 and 1.

In a similar way counterexamples can be found in the regime $d>2$ and $p>1$.
For $p\ge 1$, we find by a similar reasoning that the reversed inequality of (\ref{eq:baba}) is equivalent to
the converse of (\ref{eq:abba1}), and therefore to the converse implication
\be\label{eq:lcl2}
\Lambda^p C \Lambda^p \le C \Longleftarrow
 (C^{1/2} \Lambda C^{1/2})^{2p} \le C^{2p}.
\ee
For $d=3$ we have found counterexamples up to $p=1.25$, but no higher.
It is therefore imaginable that the second majorisation inequality in Theorem 1 could be valid under more general
conditions, e.g.\ for $p\ge2$ perhaps. For the time being, this problem is still open.
%%%%%%%%%%%%%%%%%%%%%%%%%%%%%%%%%%%%%%%%%%%%%%%%%%%%%%%%%%%%%%%%%%%%%%%%%%%%%%%%%%%%%%%%%%%%%
%------------------------------------------------------------- BIBLIOGRAPHY

%%%%%%%%%%%%%%%%%%%%%%%%%%%%%%%%%%%%%%%%%%%%%%%%%%%%%%%%%%%%%%%%%%%
\end{document}